\documentclass[10pt,a4paper]{article}

\date{\today}

\usepackage{latexsym,amssymb,amsmath}
\usepackage[latin1]{inputenc}
\usepackage[T1]{fontenc}
\usepackage[english]{babel}
\usepackage{subfig,graphicx}
\usepackage{multirow}
\usepackage{natbib}

\def\RR{\mbox{I\hskip -0.177em R}}
\def\PP{\mbox{I\hskip -0.177em P}}
\def\P{\mbox{I\hskip -0.177em P}}
\def\NN{\mbox{I\hskip -0.177em N}}

\def\E{\mbox{I\hskip -0.177em E}}
\def\II{\mbox{\rm \Large 1\hskip -0.353em 1}}

\RequirePackage{theorem}

\newtheorem{thm}{\noindent Theorem}[section]
\newtheorem{prop}{\noindent Proposition}[section]

\begin{document}

\title{A comparison of statistical models for short categorical or ordinal time series with applications in ecology}

\author{
No\"elle Bru \footnote{Universit\'e de Pau et des Pays de l'Adour, Laboratoire de Math\'ematiques Appliqu\'ees - UMR CNRS 5142 / IUT STID, Avenue de l'Universit\'e, 64013 Pau cedex, France. Email: {\tt nbru@univ-pau.fr}} \and Laurence Despres \footnote{Universit\'e Joseph Fourier, Laboratoire d'Ecologie Alpine (LECA) - UMR CNRS 5553, BP 53, 38041 Grenoble Cedex  09, France. Email: {\tt laurence.despres@ujf-grenoble.fr}} \and Christian Paroissin \footnote{Universit\'e de Pau et des Pays de l'Adour, Laboratoire de Math\'ematiques Appliqu\'ees - UMR CNRS 5142, Avenue de l'Universit\'e, 64013 Pau cedex, France. Email: {\tt cparoiss@univ-pau.fr}} }

\date{}

\maketitle

\vskip 1cm
\noindent {\it Last version:} \today

\vskip 5mm
\noindent {\it Running title:} Statistical models for categorical data time series in ecology.

\vskip 1cm
\begin{abstract}
We study two statistical models for short-length categorical (or ordinal) time series. The first one is a regression model based on generalized linear model. The second one is a parametrized Markovian model, particularizing the discrete autoregressive model to the case of categorical data. These models are used to analyze two data-sets: annual larch cone production and weekly planktonic abundance.
\end{abstract}



\newpage

\section{Introduction}

The quest to understand mechanisms behind the temporal dynamics of a natural population (animal or plant)  always yields useful information for ecological biodiversity management. The present work was motivated by the analysis of time-fluctuations of an ecological time series, the annual larch cone production. Because of the impracticability of quantitative evaluation of such population size, only semi-quantitative data are often available numerically. Data are coded with finite ordered categories or levels. From this some natural questions arise. Among them, the following one is of crucial interest here: do lagged values determine future production? The same basic problems occur as for classical quantitative time series but here the greatest difficulty stems from the nature of the studied process (as previously mentioned, working on categorical variables induces many difficulties since most of the notions used for quantitative variables have no more sense in such context).
\\[1ex]
Statistical models have been useful instruments for testing hypothesis concerning the mechanisms behind temporal evolution and to characterize temporal patterns. Two models are used throughout this paper to achieve this goal: a regression model \citep{FK1,FK2} and a parametrized Markovian model \citep{JacobsLewis1978a}. The first one is a regression model for categorical time series which is based on generalized linear regression theory \citep{McC2}. Such model extend linear models to accommodate both non-normal response distributions (which is the case in the study of categorical data) and transformations to linearity. So, applying a generalized linear models consists in two choices: a family of probability distribution and a link function between the response and the predictors. For categorical data some widely used models are: multinomial-logit \citep{Ag} and cumulative odds models \citep{McC1}. Adaptation of such models to categorical times series is easy to do putting past observations at different lags as categorical predictors of the response at time $t$ \citep{FK1,FK2}. The second one is indeed an adaption of the discrete auto-regressive (DAR) model introduced by Jacobs and Lewis \citeyearpar{JacobsLewis1978a} to the case of categorical time series. As noticed by McKenzie \citeyearpar{McKenzie}, DAR models <<~would be more suited to modelling dependent sequences of categorical observations, but this does not seem to have been attempted yet~>>. To the best of our knowledge, no advance in this direction is made since the paper of McKenzie.
\\
These two models have some advantages and some disadvantages which are not necessary the same, implying a complementarity between these two approaches. Among the common advantages, the main one is that they are easy to be interpreted by the practitioners. Since most of the time series in ecology are short-length (for a statistical purpose), we have to consider only models involving a reasonable number of parameters. That is the reason why we will focus on one order lagged model (even these models can be extended easily to large order lag values). Among the inconvenient of the DAR model, the main one is the stationarity of the time series, which can not be checked by any statistical tests (see \citep{MH} for a discussion about several notions of stationarity for categorical time series). However it allows us to derive a simple model for taking into account missing values (a contrario to the regression model, the DAR can not treat directly the case of missing values). Our approach differs highly of the one recently proposed by Bandt \citeyearpar{Bandt}. Indeed he considers a continuous-state, but non-Gaussian, time series and its analysis relies only on the ordinal property of $\RR$. Moreover his methodology requires a long time series, which is not realistic in many real cases.
\\[1ex]
Motivation of the present work is the analysis of time series of larch cone production data in spatially disjoint locations in order to determine some temporal patterns of larch cone production dynamic at different locations and to discuss some kind of spatial synchrony. Data are detailed in the first section. Next section~3 is devoted to present two regression models: one for categorical time series and one for ordinal time series. These two regression models have been studied by Fokianos and Kedem \citeyearpar{FK1,FK2}. In section 4 we adapt the one order discrete auto-regressive model to the context of categorical data (the ordinal characteristic is not taken into account in this model). In particular we develop independence tests and estimators of the various parameters of the model. In section 5, we apply these models to two real data sets: the first one deals with annual larch cone production (over 31 years) and the second one with weekly planktonic abundance (during one year). Last section is devoted to conclusion and discussion.

\section{Motivations}

\begin{figure}\centering
\subfloat[Ayes 2200]
{\includegraphics[width=7cm]{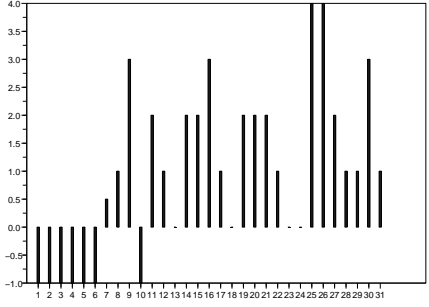}\label{fig:ayes}} 
\qquad 
\subfloat[Montgen\`evre 2200]
{\includegraphics[width=7cm]{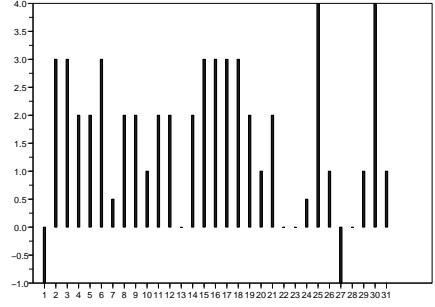}\label{fig:mntg}} 
\\
\subfloat[N\'evache 1800]
{\includegraphics[width=7cm]{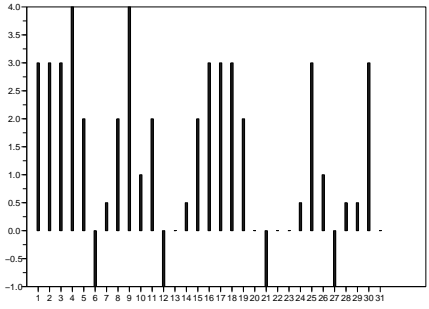}\label{fig:neva}} 
\qquad 
\subfloat[Prorel 1800]
{\includegraphics[width=7cm]{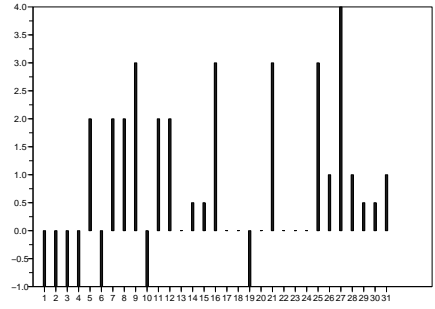}\label{fig:prorel}} 
\caption{Annual larch production in four sites in the Southern Alps}
\label{fig:larch}
\end{figure}

The masting is the intermittent synchronous production of seed crops by a plant population \citep{Kelly}. It often shows an evolved strategy related to others environmental masting patterns such as rainfall, temperatures, \ldots\, Thus variability in seed production according to past values is a good descriptor of environmental changes in climate for example. The information arising from the characterization of temporal patterns on such time series could be used to infer role of environmental parameters and other mechanisms \citep{Price}. The data accounting cone production were registered for 31 years on four valleys located in the Southern French Alps (in the same area of the Alps called "Brian\c{c}onnais"). Here we will consider four different sites selected to be at comparable altitudes (ranging from 1800m to 2200m): Ayes (altitude: 2200 meters), Montgen\`evre (altitude: 2200 meters), N\'evache (altitude: 1800 meters) and Prorel (altitude: 1800 meters). Cone production at a given site was roughly estimated at the beginning of the cone development by counting cones along one meter of branch for at least one hundred randomly selected trees. The intensity of larch cone production at any site was then classified into six classes \citep{Roques} from no cones (coded 0) to very heavy crop i.e. more than two hundred cones per tree (coded 4). Annual cone production is considered to be the realization of an ordinal time series with values $\{0, 0.5, 1, 2, 3, 4\}$ corresponding to a scale classification endowed with a natural ordering. Data are plotted on figure~\ref{fig:larch}.
\\[1ex]
When studying the dynamic of the larch population on each sampling sites, a first step could be to identify temporal patterns of cone production and then to compare each patterns from one site to others to conclude or not at a spatial synchrony on a "short" regional spatial scale \citep{Lieb}. However, the observed series in figure~\ref{fig:larch} do not exhibit obviously the presence of such patterns. The salient features of the series are:  no seasonality, high location to location variability with respect of duration and beginning of intensive larch cone production, presence of missing values, \ldots\, However visual remarks should be considered carefully. 
\\[1ex]
Such time series could appear as too short-length for the statistician who generally needs a lot of information to infer on a phenomenon but the data are collected from 1975 to 2005, which corresponds to an entire career of a biologist!

\section{Regression models for categorical and ordinal time series}

The model used here is a generalization of classical regression models to the case of time-dependent categorical observations and was studied by \citep{Kauffmann} (see also \citep{FK1,FK2} for a good summary of the main theoretical aspects).

\subsection{Introduction to generalized linear models for qualitative time series}

Assume that the observed series is a particular realization of the stochastic process in discrete time $\{Y_t\}$ which will be described below. Values of $Y_{t}$ are supposed to belong to a finite set $E=\{1, \ldots, k\}$ of $k$ ordered or not categories. Because we are interested in temporal dependence between successive observations, we condition on the observed past. For any positive integer $l$, let us denote by ${\cal{F}}_{t-l}$ the $\sigma$-field generated by $Y_{t-1}, Y_{t-2},\ldots, Y_{t-l}$. Let ${\bf Y}_t = (Y_{t,1}, \ldots, Y_{t,k-1})'$ where $Y_{t,j}$ equals to 1 if the $j$-th is observed at time $t$ and $0$ otherwise. The analysis of time series based on a generalized linear model suppose that the response variable is influenced by its past values which are viewed as predictors influencing the distribution of $Y_{t}$ by way of a transformation of a linear combination. 
\begin{equation*}
\E\left[Y_{t}|{\cal{F}}_{t-l}\right]={\bf h}\left({\bf Y}'_{t-l}\beta\right),
\end{equation*}
where $l$ is the order of the lag time and $\bf{Y}'_{t-l}$ is the covariate matrix containing the lagged values of the response variable until lag $l$. In the following we will focus on $l \in \{0,1,2\}$ (since we aim at treating short-length time series). The vector $\beta$ is a vector of time-invariant parameters to be estimate which will reflect the intensity of the dependency between the response and its past. Because the response variable is a categorical time series, we have the following relation:
\begin{equation*}
\pi_{t,j,l}= \E[Y_{t,j}|{\cal{F}}_{t-l}]= \P(Y_{t,j}=1|{\cal{F}}_{t-l}),
\end{equation*}
for every $j \in \{1,\ldots,k\}$ and every $t \{1,\ldots,n\}$, where $\pi_{t,j,l}$ is a transition probability. Let $\pi_{t,\cdot,l}=(\pi_{t,1,l}, \ldots, \pi_{t,k-1,l})$. Some adequate regression models for categorical data falls in the family of generalized linear models which links vector of transition probabilities of the response vector to the covariate process through the equation:
\begin{equation}
\pi_{t,\cdot,l} := \pi_{t,\cdot,l}(\beta)={\bf h}({\bf Y}'_{t-l}\beta) \qquad {\mbox{or}} \qquad  {\bf h}^{-1}(\pi_{t,\cdot,l}(\beta))={\bf Y}'_{t-l}\beta \;.
\end{equation}
In other words, the study of having the response $Y_{t}=j$ at time $t$ is equivalent to carry out a regression on covariates which are the lagged values of the categorical response process. This model is also called a Markov regression model for categorical time series. The function $\bf{h}$ is called the inverse link function and is related to a link function that describes how the mean depends on the linear predictors. For each response distribution there exists a variety of link functions to connect the mean with the linear predictor. The use of a generalized linear model is the choice of a combination of response distribution and a link function.

\subsection{On the choice of the link function}

The link function should adapted to the type of data \citep{FK2}:

\begin{itemize}
\item Nominal data: the most commonly used model for categorical (or nominal) data is the multinomial logit model \citep{Ag}:
\begin{equation}
\pi_{t,j,l}(\beta)=\frac{\exp(\beta'_j y_{t-l})}{1+\sum_{q=1}^{k-1}{\exp(\beta'_q y_{t-l})}} \;,
\end{equation}
for any $j \in \{1,\ldots,k-1\}$. This equation also defines log-odds ratios relative to $\pi_{tm}$ by:
\begin{equation}
\log\left(\frac{\pi_{t,j,l}}{\pi_{t,k,l}}\right)= \beta'_j {\bf y}_{t-l} \;.
\end{equation}

\item Ordinal data: since data are ordinal, its is more convenient to model the cumulative probability function of $Y_{t}$. For ordered categorical time series a reasonable choice of link function is the logistic distribution one which leads to the proportional odds model \citep{McC1}:
\begin{equation}\label{InvFunc}
{\bf h}^{-1}(x) =\frac{1}{1+\exp(-x)} \;.
\end{equation}
It follows that the link function is:
\begin{equation}\label{FuncLink}
\log\left(\frac{P\left[Y_{t}\leq j|{\cal{F}}_{t-l}\right]}{P\left[Y_{t}>j|{\cal{F}}_{t-l}\right]}\right) ={\bf Y}'_{t-l}\beta  \;.
\end{equation}
\end{itemize}

\subsection{Parameters estimation and global adequacy criteria}

Since the joint distribution of response and covariates is often not easy to establish, the likelihood methods are not applicable to estimate the vector of regression coefficients $\bf{\beta}$. As we are interested in the estimation of the effects of the covariates on the response, we can use the inference theory based on partial likelihood function. The reader can refer to \citep{FK1,Menard2} for more details and application. The partial likelihood method leads to non linear equations system.  Multinomial models were fitted using the function {\tt multinom} from library section {\tt nnet} on {\tt R}. Proportional-odds logistic regression models were fitted using the function {\tt polr} from library section {\tt MASS} on {\tt R}. The vector of parameters of this model $\bf{\beta}$ is estimated using an iterative weighted least squares {\tt IWLS} \citep{S,StatS}. 
\\[1ex]
The analysis of the global adequacy and goodness of fit of such models to the data is discussed using the Akaike's information criterion (AIC) which also allows to compare several models. The values of this criterion depends on the number of model parameters and penalizes models with large number of parameters. Such consideration is important in the study of short time series where the number of parameters can be rapidly equal to the length of the time series. The chosen model is the one which minimizes the value of AIC among the others. 
\\[1ex]
In this preliminary work, no detailed analysis of the residuals of the models is done. Such analysis is important to assess the goodness of fit between the chosen model and the observed data but was not the priority of this paper.

\section{Discrete auto-regressive model and categorical data}

The discrete auto-regressive (DAR) model introduced by Jacobs and Lewis \citeyearpar{JacobsLewis1978a,JacobsLewis1978b,JacobsLewis1978c} is used here to model categorical data. Some independence tests are developed, using either the Markov property or runs properties. Estimators of the parameters are studied in the precise context of categorical data. Simulated data are used to illustrate numerically the quality of these estimators.

\subsection{Introduction and model}

In a series of papers, Jacobs and Lewis \citeyearpar{JacobsLewis1978a,JacobsLewis1978b,JacobsLewis1978c} introduced and studied time series models for discrete variables. Among them we will focus here on the discrete auto-regressive of order 1, denoted by DAR(1). Such process $\{X_t\}$ is a discrete-time stochastic process with values on a finite ordered set $E= \{1, \ldots, k\}$ and is defined as follows:
\begin{equation*}
\forall t>0\;, \quad X_t = V_t X_{t-1} + (1-V_t) Z_t \;,
\end{equation*}
where $\{V_t\}$ is a sequence of iid Bernoulli random variables with parameter $\alpha \in [0;1]$ and $\{Z_t\}$ is a sequence of iid random variables having the distribution $\pi$ on $E$, the two sequences being independent. Moreover we will assume that $X_0$ is distributed according to the distribution $\pi$, implying that the process $\{X_t\}$ is stationary. The case of $\alpha=1$ is not interesting since $X_t=X_0$, with probability 1, for any $t$. The case of $\alpha=0$ means that the process $\{X_t\}$ is simply a sequence of iid random variables having distribution $\pi$. Hence the parameter $\alpha$ could be interpreted as follows: the nearest to 0 $\alpha$ is, the more <<~independent~>> the sequence $\{X_t\}$ is. Indeed, for all $h \in \NN$, $\rho(h)= \alpha^h$ is the auto-correlation function of a DAR(1) process. Hence DAR(1) models can be used to describe a situation of short range dependency with high correlation. It is easy to prove that stochastic process $\{X_t\}$ as defined above is a Markov chains on $E$ with transition matrix $P$ given by the following equation:
\begin{equation*}
P= \alpha I + (1-\alpha)Q \;,
\end{equation*}
where ${}^t Q = [{}^t\pi|\cdots|{}^t\pi]$. Such Markov chain admits obviously a unique stationary probability distribution which is $\pi$. One can easily deduce the $h$-th power of $Q$ and $P$: for all $h \geqslant 1$, $Q^h= Q$ and $P^h= \alpha^h I + (1-\alpha^h)Q$, illustrating one more times the role of $\alpha$.
\\[1ex]
This stochastic process could be generalized to higher order leading to the DAR($p$) model. In fact, these models appear themselves to be a special case of mixture transition distribution (MTD) model introduced by Raftery \citeyearpar{R}. Thus DAR($p$) can be viewed as an alternative to MTD model. According to Raftery, a MTD model fits better data in general than a DAR($p$) one, especially for $p \geqslant 3$. However here we will prefer to use a DAR(1) model since it has the following advantages over the MTD model: 1) the two parameters $\alpha$ and $\pi$ play different roles: $\alpha$ is related to the correlation whereas $\pi$ is the stationary distribution; 2) these models involve generally a reasonable number of parameters (more parsimonious) especially when few data are available; 3) parameters could be easily interpreted by a practitioner. But the special case of DAR(1) model presents the disadvantage of being restrictive over the transition matrix. 
\\[1ex]
Here we are interested on the use of such stochastic processes for modeling categorical variables (here the $k$ different modalities are encoded by using the $k$ first positive integers). It implies that many characteristics of these processes have no sense in such context, as it is the case for the auto-correlation function (see above). Thus estimators developed by Jacobs and Lewis (\citeyear[see pages 28--30]{JacobsLewis1983}) cannot be used. Hence we address the statistical problem of estimating the parameters in a DAR(1) model in presence of categorical data. Assume we observe $X_0, \ldots, X_n$ for a fixed value $n>0$. First we will test whether $\{X_t\}$ is a sequence of independent random variables ($\alpha=0$) or not. In a second step we will estimate all the parameters of the model: $\alpha$ and $\pi$. Then we propose a very simple model in order to consider the case of missing observations.

\subsection{Independence tests}

In this section we aim at testing whether $\{X_t\}$ is a sequence of independent random variables ($\alpha=0$) or not. Two ways will be investigated. The first one will use the Markov property of the DAR(1) model and the second one will be based on runs property. Anyway, all along this section, the null hypothesis $H_0$ will be <<~$\alpha=0$~>> and the alternative hypothesis will be  <<~$\alpha \neq 0$~>>.

\paragraph{$\chi^2$ test based on the Markov property} The following test is a classical test for Markov chain (see \citep{RSW} for an illustration in DNA analysis context). We only use the fact that $\{X_t\}$ is a Markov chain, but not the particular structure of its transition matrix. Classical results on Markov chain inference leads to the following estimate for the transition matrix $P$:
\begin{equation*}
\widehat{P}_{j,j'}= \frac{N^n_{j,j'}}{N^n_{j,\cdot}} \;,
\end{equation*}
where $N^n_{j,j'} = \sum_{i=1}^{n} \II_{\{X_{i-1}=j, X_{i}=j'\}}$ and $N^n_{j,\cdot} = \sum_{j' \in E} N^n_{j,j'} = \sum_{i=1}^{n} \II_{\{X_{i-1}=j\}}$. In other words $N^n_{j,j'}$ is the number of jumps from state $j$ to state $j'$ and $N^n_{j,\cdot}$ is the number of visits of state $j$, in the sequence of observations $X_0, \ldots, X_n$.
\\[1ex]
The null hypothesis can rephrased as follows: $P_{j,j'}=P_{j,\cdot}P_{\cdot,j'}$, for any $(j,j') \in E^2$. Under $H_0$, the maximum-likelihood estimate of $P_{j,j'}$ is:
\begin{equation*}
\widehat{P}_{j,j'}= \widehat{P}_{j,\cdot}\widehat{P}_{\cdot,j'}=
\frac{N^n_{j,\cdot}}{n-1}\frac{N^n_{\cdot,j'}}{n-1} \;,
\end{equation*}
where $N^n_{\cdot,j'} = \sum_{j \in E} N^n_{j,j'} = \sum_{i=1}^{n} \II_{\{X_{i}=j'\}}$. Hence one has to consider the following statistics $C^2$:
\begin{equation*}
C^2 = \sum_{j\in E}\sum_{j'\in E} \frac{[N^n_{j,j'}-N^n_{j,\cdot}N^n_{\cdot,j'}/(n-1)]^2}{N^n_{j,\cdot}N^n_{\cdot,j'}/(n-1)} \;.
\end{equation*}

\begin{thm}
Under the null hypothesis, $C^2  \xrightarrow[n \rightarrow \infty]{d} \chi^2_{(k-1)^2}$ \;.
\end{thm}

Some well-known practical restrictions exist in order to be able to apply this test. As example, one can require that $\widehat{P}_{j,j'}>5\%$, for any $(j,j') \in E^2$.

\paragraph{Tests based on runs property} Unfortunately we cannot compute the power of the previous test, that is the reason we will now consider a second family of tests. These tests will be based on runs property of the model. Runs in sequence of iid Bernoulli distributions are studied for a very long time: this problem seems to be considered for the first time by Abraham de Moivre in 1756 (problem {\sc LXXIV} in his book {\em The Doctrine of Chances}). For an historical perspective, see the introduction of the Part I of \citep{Mood}. Most of the existing papers deal with the case of Bernoulli random variables, but here we are indeed interested in the general discrete case. Few extensions were made in this direction. To the best of our knowledge, Mood \citeyearpar{Mood} is the first one who studied it.
\\[2ex]
A run can be defined as follows: it is a consecutive sub-sequence of identical values in a sequence of random numbers.  For any $j \in E$, let us denote by $R_{j,n}^i$ the number of non-overlapping $j$-runs of length $i$ in the sequence $X_1, \ldots, X_n$:
\begin{equation*}
R_{j,n}^i = |\{ m \;;\; X_{m-1} \neq j \,,,\, X_{m}= j \,,,\, \ldots \,,,\, X_{m+i-1}= j \,,,\, X_{m+i} \neq j \}| \;.
\end{equation*}
Let us now define the number $R_{j,n}$ of $j$-runs and the total number $R_{n}$ of non-overlapping runs:
\begin{equation*}
R_{j,n} = \sum_{i=1}^{n} R_{j,n}^i \;,
\qquad
{\mbox{and}}
\qquad
R_{n} = \sum_{j \in E} R_{j,n} \;.
\end{equation*}

Mood \citeyearpar{Mood} obtained the limiting distribution of $R_{n}$ after renormalization. Two cases have be distinguished: $k=2$ and $k>2$.

\begin{thm}
(corollary~5 p.~390 and corollary~3 p.~392 in \citep{Mood})
\begin{enumerate}
\item If $k=2$, $\displaystyle{\frac{R_{n}- 2n\pi_1\pi_2}{2\sqrt{n\pi_1\pi_2(1-3\pi_1\pi_2)}}  \xrightarrow[n \rightarrow \infty]{d} {\cal{N}}(0,1)}$.
\item If $k>2$, $\frac{1}{\sqrt{n}}\left( R_{n}- n\left(1-\sum_{j\in E}\pi_j^2\right) \right) \xrightarrow[n \rightarrow \infty]{d} {\cal{N}}(0,\sigma^2_{\pi})$, with $\sigma^2_{\pi} = \sum_{j \in E} \pi_j^2 + 2\sum_{j \in E} \pi_j^3 - 3( \sum_{j \in E} \pi_j^2 )^2$.
\end{enumerate}
\end{thm}

One can check that in both cases the asymptotic variance is degenerated if and only if there exists $j \in E$ such that $\pi_j=1$ (and then, for all $j' \neq j$, $\pi_{j'}=0$), then the variance is degenerated. These convergence results could be used to construct an asymptotic test. 
\\[2ex]
An alternative solution could be to consider the longest run in the sequence $X_1, \ldots, X_n$. Indeed there exists many works dealing with the case of either independent trials or Markovian trials. Let us denote by $L_n$ the longest length of all runs in the sequence $X_1, \ldots, X_n$. Using the previous notations, we have:
\begin{equation*}
L_n = \max \{ i \;;\; \exists j \; s.t. \; R_{j,n}^i>0 \} \;.
\end{equation*}
Vaggelatou \citeyearpar{Vaggelatou} studied this random variable in the case of multi-state Markovian trials. It requires that $\{X_t\}$ is an irreducible and aperiodic Markov chain on a finite state space $E$ with transition probability matrix $P$ and unique stationary measure $\pi$. The Markov chain induced by a DAR(1) model is irreducible and aperiodic if $\pi>0$ (meaning that all the components of $\pi$ are strictly positive) and $\alpha \neq 1$ (see chapter {\sc xv} of \citep{Feller1}). Let us define the two following quantities:
\begin{equation*}
\rho = \max_{j \in E} P_{jj}
\quad {\mbox{and}} \quad
\pi_{\rho} = \sum_{j \in E \::\: P_{jj}= \rho}\pi_j \;.
\end{equation*}
If $\rho<1$, then Vaggelatou proved the following asymptotic result (theorem~1 in \citep{Vaggelatou}):
\begin{equation}\label{eqn:V}
\PP( L_n - [\log_{1/\rho} n] < x ) = \exp\left\{-n (1-\rho)\pi_{\rho} \rho^{[\log_{1/\rho} n]+x-1}\right\} + o(1) \;,
\end{equation}
where $[\cdot]$ denotes the integer part and $o(1)$ means that the residual term is <<~small~>> in regard with $n$. This result extends the classical one obtained many years ago by Gon\u{c}arov \citeyearpar{G} in the case of iid Bernoulli trials. In both case, $L_n - [\log_{1/\rho} n]$ does not have a limit distribution, but only certain sub-sequence; for instance, theorem~2 in \citep{Vaggelatou} gives a case where the sub-sequence converges in distribution to the Gumbel distribution. We will use theorem~1 (and not theorem~2) to construct a third (and last) test since we have not enough observations in real situation. Let us denote by $\rho_0$ and $\rho_1$ the value of $\rho$ respectively under the null and the alternative hypothesis. Under the null hypothesis, $P$ will be equal to the matrix $Q$ as defined in the introduction: $\forall (j,j') \in E^2$, $P_{jj'} = \pi_{j'}$ (the transition probability from state $j$ to state $j'$ does not depend on $j$). So we have that $\rho_0 = \max\{\pi_j \;;\; j \in E\}$: $\rho_0<1$ if and only if, for any $j \in E$, $\pi_j<1$. Under the alternative hypothesis, $\rho_1 = \max\{P_{jj} \;;\; j \in E\} =  \max\{\alpha+ (1-\alpha)\pi_j \;;\; j \in E\}$: $\rho_1<1$ if and only if $\alpha<1$ (it is initially assumed) and  for any $j \in E$, $\pi_j<1$. Thus we find the same condition in both cases and this assumption is the same as for the previous test. From now we will assume that $\pi>0$ in addition to the previous assumptions (let us recall that we already assume that $\alpha \neq 1$). Using theorem~1 of Vaggelatou \citeyearpar{Vaggelatou}, one could obtain an asymptotic confidence interval with a prescribed confidence level $\varepsilon \in (0;1)$:
\begin{equation*}
\PP_{H_0} \left( \tilde{L}_n \in \bar{W}_{\varepsilon,n} = \left[ \log_{\rho_0} \left( - \frac{\ln (\varepsilon/2)}{n(1-\rho_0)\pi_{\rho_0}}\right) \;;\;  \log_{\rho_0} \left( - \frac{\ln (1-\varepsilon/2)}{n(1-\rho_0)\pi_{\rho_0}}\right) \right] \right) = 1-\varepsilon \;,
\end{equation*}
$\tilde{L}_n = L_n - 1$ (notice that it is corresponding sometimes to the definition of runs: see for instance \citep{JacobsLewis1978a}). It follows that the power $\Pi_{\varepsilon}$ of this test is:
\begin{equation*}
\Pi_{\varepsilon} = 1+\PP_{H_1} \left(\tilde{L}_n <\log_{\rho_0} \left( - \frac{\ln (\varepsilon/2)}{n(1-\rho_0)\pi_{\rho_0}}\right)  \right)  - \PP_{H_1} \left(\tilde{L}_n <\log_{\rho_0} \left( - \frac{\ln (1-\varepsilon/2)}{n(1-\rho_0)\pi_{\rho_0}}\right)  \right) \;.
\end{equation*}
To compute $\Pi_{\varepsilon}$, one has to use equation~(\ref{eqn:V}) above.

\subsection{Parameters estimations}

The two parameters $\pi$ and $\alpha$ of such DAR(1) model could be estimated separately since by construction they play different role.

\paragraph{Estimations of $\pi$} For any $j \in E$ and for any $i \in \{1, \ldots, n\}$, let $Z_{ij} = \II_{\{X_i=j\}}$: these random variables have the Bernoulli distribution with parameter $\pi_j$. A natural unbiased estimator of $\pi_j$ is therefore:
\begin{equation*}
\widehat{\pi}_j = \frac{1}{n} \sum_{i=1}^n Z_{ij} \;.
\end{equation*}
Moreover, using the expression of $P^h$ given in the introduction, one can easily derive the variance of $\widehat{\pi}_j$:
\begin{equation*}
{\mbox{Var}}[\widehat{\pi}_j] = \frac{1}{n} \pi_j (1-\pi_j) + \frac{2}{n^2} (1-\pi_j)\pi_j V_n(\alpha) \;,
\end{equation*}
and the covariance between $\widehat{\pi}_{j}$ and $\widehat{\pi}_{j'}$ (with $j'\neq j$):
\begin{equation*}
{\mbox{Cov}}[\widehat{\pi}_j,\widehat{\pi}_{j'}] = -\frac{2}{n^2} \pi_{j'} \pi_j V_n(\alpha)  \;,
\end{equation*}
with $V_n(\alpha) = \sum_{h=1}^{n} (n-h)\alpha^h$. Applying formula (0.113) in \citep{GR} (arithmetico-geometric progression), we obtain the following  expression for $V_n(\alpha)$:
\begin{equation*}
V_n(\alpha) = \frac{n-\alpha^n}{1-\alpha} - \frac{\alpha(1-\alpha^{n-1})}{(1-\alpha)^2} - n \;.
\end{equation*}
This leads to the following limit for the variances and the covariances:
\begin{equation*}
\lim_{n \rightarrow \infty}n{\mbox{Var}}[\widehat{\pi}_j] =\frac{1+\alpha}{1-\alpha}  \pi_j (1-\pi_j) \;,
\end{equation*}
and:
\begin{equation*}
\lim_{n \rightarrow \infty}n  {\mbox{Cov}}[\widehat{\pi}_j,\widehat{\pi}_{j'}] = -\frac{2\alpha}{1-\alpha} \pi_{j'} \pi_j   \;.
\end{equation*}
As a consequence of Bienaym\'e-Chebychev inequality, the asymptotic result on the variance implies that $\widehat{\pi}_j$ is consistent:
\begin{prop}
For any $\alpha \in [0,1)$, $\widehat{\pi}_j \xrightarrow[n \rightarrow \infty]{Pr} \pi_j$.
\end{prop}

Moreover one can prove the following central limit theorem for $\widehat{\pi}_j$ by application of the ergodic theorem for Markov chain \citep{Jones} and Slutsky theorem:

\begin{thm}
For any $\alpha \in [0,1)$, $\displaystyle{\sqrt{n} \frac{\widehat{\pi}_j - \pi_j}{\sqrt{\widehat{\pi}_j(1-\widehat{\pi}_j)}} \xrightarrow[n \rightarrow \infty]{d} {\cal{N}}\left(0,\frac{1+\alpha}{1-\alpha}\right)}$.
\end{thm}

When $\alpha=0$, the asymptotic variance equals to 1 as it is well-known for Bernoulli trials. The largest $\alpha$ is, the larger the asymptotic variance is. Since the asymptotic depends on $\alpha$ which is generally unknown, we cannot yet use the last proposition to construct confidence interval. In order to do it, we will need an consistent estimator of $\alpha$.

\paragraph{Estimations of $\alpha$} We will first consider the maximum likelihood estimator of $\alpha$, assuming that $\pi$ is known. Since $\{X_t\}$ is a Markov chain, the log-likelihood is:
\begin{equation*}
{\cal{L}}(X_1, \ldots, X_n;\alpha)  = \sum_{(j,j') \in E^2} N_{j,j'} \log P_{j,j'}(\alpha') \;,
\end{equation*}
where $N_{j,j'}$ is defined as in section~2 and where $P_{j,j'}(\alpha)$ is the transition probability from state $j$ to state $j'$ (which only depends on $\alpha$). Replacing the expression of transition probabilities, we obtain that:
\begin{equation*}
{\cal{L}}(X_1, \ldots, X_n;\alpha)  = \sum_{j \in E} \left( N_{j,j}  \log (\alpha + (1-\alpha)\pi_{j}) + \sum_{j' \in E \setminus \{j\}} N_{j,j'}\log ((1-\alpha)\pi_{j'}) \right) \;.
\end{equation*}
It follows that the maximum likelihood estimator $\alpha_1^*$ of $\alpha$ is the solution of the following equation:
\begin{equation*}
\frac{1}{n} \sum_{j \in E} \frac{N^n_{j,j}}{\alpha+(1-\alpha)\pi_j} = 1 \;.
\end{equation*}
When $\pi$ is unknown, one can use the estimation given above and so the plug-in estimator $\widehat{\alpha_1}$ of $\alpha$ is the solution of the following equation:
\begin{equation*}
\frac{1}{n} \sum_{j \in E} \frac{N^n_{j,j}}{\alpha+(1-\alpha)\widehat{\pi}_j} = 1 \;.
\end{equation*}
Unfortunately we cannot derive an explicit expression of $\widehat{\alpha_1}$. An alternate possible way is to minimize the following function:
\begin{equation*}
Q_\alpha = \sum_{(j,j') \in E^2} (\widehat{P}_{j,j'} - P_{j,j'}(\alpha))^2 \;.
\end{equation*}
Solving this optimization problem leads to an explicit expression ${\alpha_2}^*$:
\begin{equation*}
\alpha_2^* = \frac{\displaystyle{\sum_{j \in E} (1-\pi_j)(\widehat{P}_{jj}-\pi_j) - \sum_{j \in E} \sum_{j' \in E \setminus \{j\}} \pi_{j'}(\widehat{P}_{jj'}-\pi_{j'})}}{\displaystyle{(k-1)\sum_{j \in E} \pi_j^2 + \sum_{j \in E}} (1-\pi_j)^2} \;.
\end{equation*}
It seems to be difficult to establish properties of this intuitive estimator. Hence it is not recommended to use it as an estimator of $\alpha$. However it could provide a possible initialization for an optimization procedure to obtain a numerical value of $\widehat{\alpha_1}$. The corresponding plug-in estimator $\widehat{\alpha_1}$ of $\alpha$ is given by the following expression:
\begin{equation*}
\widehat{\alpha_2} = \frac{\displaystyle{\sum_{j \in E} (1-\widehat{\pi}_j)(\widehat{P}_{jj}-\widehat{\pi}_j) - \sum_{j \in E} \sum_{j' \in E \setminus \{j\}} \widehat{\pi}_{j'}(\widehat{P}_{jj'}-\widehat{\pi}_{j'})}}{\displaystyle{(k-1)\sum_{j \in E} \widehat{\pi}_j^2 + \sum_{j \in E}} (1-\widehat{\pi}_j)^2} \;.
\end{equation*}

\paragraph{Simulated data} The estimators developed above are applied on simulated data in order to evaluate numerically their performance. We simulate data with various values of $\alpha$, $\pi$ and $n$:
\begin{itemize}
\item  $\pi=(\tfrac{1}{2},\tfrac{1}{2})$, $\pi=(\tfrac{1}{3},\tfrac{2}{3})$, $\pi=(\tfrac{1}{3},\tfrac{1}{3},\tfrac{1}{3})$ and $\pi=(\tfrac{1}{4},\tfrac{1}{2},\tfrac{1}{4})$
\item $\alpha \in \{0.1;0.2;0.5;0.8;0.9\}$.
\item $n \in \{50;100;500\}$.
\end{itemize}
Hence two cases are studied: $k=2$ and $k=3$. For each values of these parameters, we simulate $m=100$ independent DAR(1) Markov chain and then we compute the estimators. Unfortunately these we have no guarantee that the two estimators of $\alpha$ belong to the unit interval. Hence we precise for each cases the number of samples on which the computations were done (as it is reasonable, the number of samples is increasing with the number $n$ of observations). Results are given in tables~\ref{tab:sim1} to~\ref{tab:sim4} (only the estimators are computed for simulated data).

\subsection{A variant with missing observations}

Sometimes categorical time series may contain some missing values/observations. Here we now propose a very simple adaptation of the DAR model in order to taking into account the missing values. Since the DAR model is stationary, it will be easy to derive similar expressions as in the initial model.
\\[1ex]
Assume that at each unit of time, the probability of a missing value equals to $\beta$ (which does not depend on the time). Hence, if we denote by $Z_t$ the values at time $t$, we have :
\begin{equation*}
Z_t = \left\{
\begin{array}{ll}
X_t & {\mbox{w.p.}} \quad   1-\beta  \\
-1 & {\mbox{w.p.}} \quad  \beta
\end{array} \right. \;,
\end{equation*}
where $-1$ is the value corresponding to a missing value and $\{X_t\}$ is DAR(1) stochastic process as described previously. Hence $\beta$ is the probability that a value is not observed : this probability is assumed to be not depending on $t$. Since $(X_t)$ is a stationary stochastic process, it follows that $(Z_t)$ is still a Markov chain, but taking values on the set $\tilde{E} =\{-1\} \cup E$. Its transition probabilities matrix $\tilde{P}$ can be expressed in function of the transition probabilities matrix $P$ of $(X_t)$:
\begin{equation*}
\tilde{P} = \left[
\begin{array}{cc}
1-\beta & \beta {}^t \pi \\
(1-\beta){\bf 1}_k & \beta P
\end{array}
\right] \;,
\end{equation*}
where ${\bf 1}_k$ is the unit vector of $\RR^k$. There is now three parameters to be estimated. Indeed $\pi$ and $\alpha$ (with the maximum likelihood method) can be estimated as previously. The extra parameter $\beta$ can be simply estimated as follows:
\begin{equation*}
\widehat{\beta} = \frac{1}{n} \sum_{i=1}^n \II_{\{X_i=-1\}} \;.
\end{equation*}

\section{Applications to ecological data}

Finally tests and estimators are applied in this section to real data. We apply the two models described previously to two real data sets. The first one deals with larch cone production (see section~2) and the second one to planktonic abundance. 

\subsection{Larch cone production}

The total number of observations was $n=31$, with $k=6$ categories. The goal of this work is to study masting on such trees. For a full description of the data, see section~2.
\\[1ex]
First we apply regression models. Tables~\ref{tab:larch-reg1} and ~\ref{tab:larch-reg2} contain values of AIC respectively for the categorical time series regression model and the ordinal one. We also indicate the number of parameters to estimate and the number of observations (these time series may contain missing values). For the first case, the independence assumption leads to the better for all sites. For the second case, model with a lag order 1 and 2 fits better for the site Ayes 2200 and Montgen\`evre 2200 while the model with only a lag of order 1 fits better for N\'evache 1800 and the independent model fits better for Prorel 1800. Comparing values of AIC, the model for ordinal time series seems to be more accurate for these data sets.
\\[1ex]
Second we apply the DAR model. Table~\ref{tab:larch-dar} contains the estimations of the three parameters for each of the four data sets. In all cases, the independence hypothesis is rejected with the two first tests, while the third one leads to accept the assumption of independence (all with a first type error at 5\%). However the power of this last test is more or less weak in all cases (ranging from above 26\% to 51\%). Thus it is reasonable to reject the assumption of independent observations.
\\[1ex]
In all cases the categorical time series regression model leads to accept the temporal independence between observations. It may be due to the fact that the parameter $\alpha$ in the DAR model is closed to zero. However this parameter can be assumed to be significantly different of zero, according to the performed tests. It is in concordance with the fact that observations are time-dependent when using the ordinal time series regression model. Time series studied here are very short-length and it may the cause that the conclusions based on the ordinal time series regression model and the ones based on the DAR models. Since few data are available, one should prefer to use the DAR models (because it involves less parameters than the other models).

\subsection{Planktonic abundance}

\begin{figure}\centering
\subfloat[Year 1987]
{\includegraphics[width=7cm]{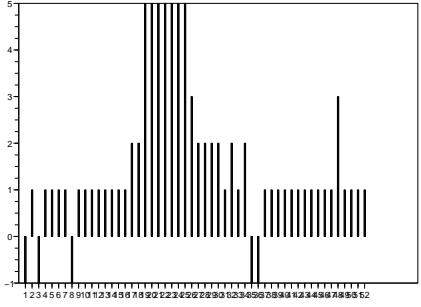}\label{fig:m1987}} 
\qquad 
\subfloat[Year 1988]
{\includegraphics[width=7cm]{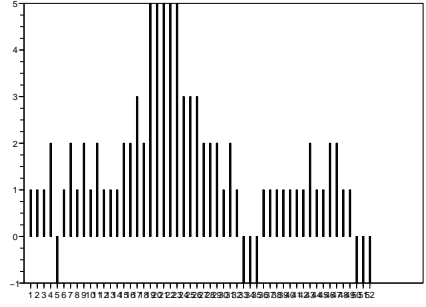}\label{fig:m1988}} 
\\
\subfloat[Year 1989]
{\includegraphics[width=7cm]{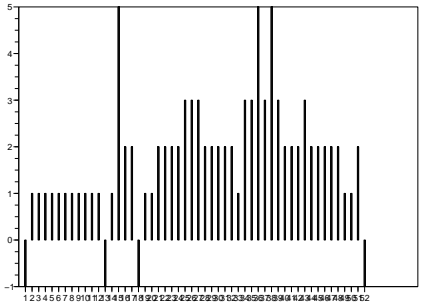}\label{fig:m1989}} 
\qquad 
\subfloat[Year 1990]
{\includegraphics[width=7cm]{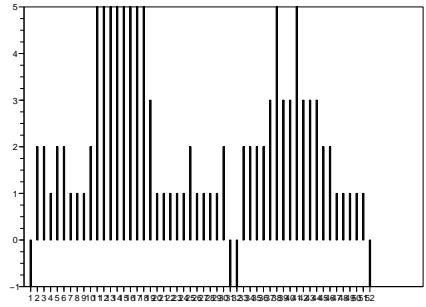}\label{fig:m1990}} 
\caption{Weekly planktonic abundance for four years}
\label{fig:menard}
\end{figure}

We now consider weekly planktonic (\emph{Thalia democratica}) abundance data. Data were kindly given by F.~M\'enard. In such context, the objective is to test and to compare the temporal patterns from one year to another. Hence we apply the two models described above for four years (1987 to 1990). It follows that each data-set is made of $n=52$ observations. Abundances were determined semi-quantitatively according to classes defined on scale of 5 values, $E=\{1,2,3,4,5\}$. The observed series are shown in figure~\ref{fig:menard} and exhibit the same problems as the larch cone production ones. Notice that the fifth category were not observed for any year. For a complete description of the data, the reader could reefer to \citep{Menard1} (see also \citep{Menard2}).
\\[1ex]
First we apply regression models. Tables~\ref{tab:menard-reg1} and~\ref{tab:menard-reg2} contain the number of parameters to be estimated, the values of AIC and the number of observations, respectively for the categorical time series regression model and the ordinal one. With the model for categorical time series, the model with one order lag fits better all the four years. With the model for categorical time series,  model with two order lag fits better for the year 1987 and 1989 while model with only a one order lag  fits better for the two other years, 1988 and 1990.
\\[1ex]
Second we apply the DAR model. Table~\ref{tab:menard-dar} contains the estimations of the three parameters for each of the four data sets. The two first tests leads to reject the null hypothesis, i.e. to reject the independence of the observations. The third test leads also to reject the independence assumption for the two last year (1989 and 1990) while the null hypothesis is accepted according to this last test for the years 1987 and 1988 (respectively with a power equal to 44\% and 66\%). Thus it is also reasonabke to reject the assumption of independent observations.
\\[1ex]
We can also analyze these data sets as one unique time series (notice that it was not possible for the previous data-set). Hence the sample size is now $n=208$. All the tests lead to reject the assumption of independent observations (with a power equal to 69\% for the last one). The last line of table~\ref{tab:menard-dar} contains the estimations of the three parameters for the whole period.
\\[1ex]
Here the situation is totally different than previously. Indeed the categorical time series regression model and the DAR model lead to the same conclusion, i.e. a one order lag dependence in the time series. One can notice that for these data sets the parameter $\alpha$ (of the DAR model) is now between 0.308 and 0.540. However the ordinal time series regression model leads in some case to a two order lag dependence. Since the number of observations is almost the twice than for the first data sets, one should rather prefer to use the ordinal time series regression model.

\section{Conclusions and discussions}

Applications to real data achieve to convince that these two complementary models are relevant for practical purpose. Based on the result obtained over real data, one can conclude that either the ordinal time-series regression model or the DAR model should be used to treat such data. Indeed in all cases the categorical time-series regression model seems not to present advantages over the two other models. The choice between the ordinal time-series regression model and the DAR model depends highly on the context, i.e. essentially on the number of observations and the number of parameters to be estimated. Fokianos and Kedem \citeyearpar{FK2} claimed that <<~the regression methodology can discover dependencies in the DNA sequence data which cannot assessed by a Markov model~>>. However data treated here can serve as a counter-example of this sentence. For the two data sets studied here, conclusions based on the regression models and the one based on the Markov model are almost identical.
\\[1ex]
From this work, one can conclude that in any case the DAR model has only few parameters to be estimated, but with an equal to number of unknown parameters (as in the example of larch cone production) one has to prefer the ordinal time-series regression model.
\\[1ex]
However both suffers of relying on assumptions or simplifications. Hence these models could be extended in the following ways:

\begin{itemize}
\item {\em Stationarity:} the DAR model is strongly stationary (in the sense defined by McGee and Harris in \citeyearpar{MH}). This assumption should be checked with any statistical tests. However no test of stationarity of a categorical time series has been developed to the best of our knowledge. Anyway when dealing with short-length time series stationary assumption is not really restrictive. Otherwise a solution could be in applying the de-trend algorithm suggested by  McGee and Harris in \citeyearpar{MH}. However their algorithm is more and more computationally complex as the number of states is increasing (in fact they mainly consider the binary case). We do not focus here on the study of the possible stationarity of a categorical time series which will be done in a future work. A major advantage of regression model is that it is not necessary to have stationarity. 

\item {\em Higher dependence order and number of parameters:} since we consider the case of short-length data, we limit our study to one or two order lagged models. Indeed both models could be applied to $p$-th order lagged models. However, for the regression model, a large value of $p$ implies a large number of parameters to estimate, that may induce some numerical instability (due to correlation between the regressors). The number of parameters in a DAR($p$) model is lower, but if $p>1$ we have no more the Markov property. 

\item {\em Environmental factors:} these models do not include environmental covariates. The regression model could easily integrate such situations, as shown by Fokianos and Kedem \citeyearpar{FK1,FK2}. For the DAR model, it is not so easy. A solution could be to consider inhomogeneous Markov chain or a state-space model.
\end{itemize}


\vspace{0.5cm}

\paragraph{Acknowledgments} 
We wish to thank Patricia Jacobs (Naval Postgraduate School, Monterey, California, USA), Ian R. Harris (Southern Methodist University, Texas, USA), Alain Latour (LabSAD, Grenoble, France), Monnie McGee (Southern Methodist University, Texas, USA) and Fr\'ed\'eric M\'enard (IRD, Montpellier, France). 
Thanks also to the numerous scientists and crew members who conducted the experimental campaigns because our analysis is based on their hard work. 
This work is a contribution to the understanding of larch cone production into an IFB (Institut Fran\c{c}ais de la Biodiversit\'e) project.


\bibliographystyle{natbib}

\bibliography{mlz}

\section*{Biographical sketches}

No\"elle Bru is an Assistant Professor at the IUT STID of Universit\'e de Pau et des Pays de l'Adour and a Researcher at the Laboratoire de Math\'ematiques Appliqu\'es de Pau. Since her PhD thesis, her topics of interest is both theoretical and applied statistics with emphasis to environmental data analysis.
\\[1ex]
Laurence Despres is an Assistant Professor at the Universit\'e Joseph Fourier and a Researcher at the Laboratoire d'Ecologie Alpine. Her main research interests are in the evolutionary ecology of species interactions and coevolution.
\\[1ex]
Christian Paroissin is an Assistant Professor at the D\'epartement Sciences et Techniques of Universit\'e de Pau et des Pays de l'Adour and a Researcher at the Laboratoire de Math\'ematiques Appliqu\'es de Pau. His topics of interest is applied probability and statistics with interest to applied contexts (engineering, theoretical computer science, biology, \ldots).


\newpage

\newpage

\newpage

\begin{table}
\begin{center}
\begin{tabular}{|c|c||c||c|c||c|c|}
\hline
$\alpha$ & $n$ & $\widehat{\pi}$ & $\widehat{\alpha_1}$ & $m_1$ & $\widehat{\alpha_2}$ & $m_2$ \\
\hline \hline
\multirow{3}{0.5cm}{0.1} & 50 & (0.509;0.491) & 0.140 & 65 & 0.149 & 72 \\
\cline{2-7}
 & 100 & (0.499;0.501) & 0.114 & 82 & 0.121 & 85 \\
\cline{2-7}
 & 500 & (0.503;0.497) & 0.092 & 100 & 0.094 & 100 \\
\hline \hline
\multirow{3}{0.5cm}{0.2} & 50 & (0.511;0.489) & 0.174 & 83 & 0.183 & 90 \\
\cline{2-7}
 & 100 & (0.510;0.490) & 0.180 & 97 & 0.192 & 97 \\
\cline{2-7}
 & 500 & (0.499;0.501) & 0.191 & 100 & 0.194 & 100 \\
\hline \hline
\multirow{3}{0.5cm}{0.5} & 50 & (0.501;0.499) & 0.451 & 99 & 0.481 & 99 \\
\cline{2-7}
 & 100 & (0.510;0.490) & 0.465 & 100 & 0.479 & 100 \\
\cline{2-7}
 & 500 & (0.499;0.501) & 0.494 & 100 & 0.497 & 100 \\
\hline \hline
\multirow{3}{0.5cm}{0.8} & 50 & (0.543;0.457) & 0.711 & 99 & 0.753 & 99 \\
\cline{2-7}
 & 100 & (0.534;0.466) & 0.765 & 100 & 0.784 & 100 \\
\cline{2-7}
 & 500 & (0.502;0.498) & 0.795 & 100 & 0.799 & 100 \\
\hline \hline
\multirow{3}{0.5cm}{0.9} & 50 & (0.585;0.415) & 0.754 & 99 & 0.798 & 99 \\
\cline{2-7}
 & 100 & (0.539;0.461) & 0.857 & 100 & 0.882 & 100 \\
\cline{2-7}
 & 500 & (0.503;0.494) & 0.890 & 100 & 0.894 & 100 \\
\hline
\end{tabular}
\caption{Results obtained with $\pi=(\tfrac{1}{2},\tfrac{1}{2})$}
\label{tab:sim1}
\end{center}
\end{table}

\begin{table}
\begin{center}
\begin{tabular}{|c|c||c||c|c||c|c|}
\hline
$\alpha$ & $n$ & $\widehat{\pi}$ & $\widehat{\alpha_1}$ & $m_1$ & $\widehat{\alpha_2}$ & $m_2$ \\
\hline \hline
\multirow{3}{0.5cm}{0.1} & 50 & (0.328;0.672) & 0.145 & 65 & 0.160 & 69 \\
\cline{2-7}
 & 100 & (0.337;0.663) & 0.120 & 74 & 0.131 & 15 \\
\cline{2-7}
 & 500 & (0.337;0.663) & 0.095 & 97 & 0.097 & 97 \\
\hline \hline
\multirow{3}{0.5cm}{0.2} & 50 & (0.335;0.665) & 0.202 & 91 & 0.226 & 91 \\
\cline{2-7}
 & 100 & (0.335;0.665) & 0.205 & 96 & 0.217 & 96 \\
\cline{2-7}
 & 500 & (0.336;0.664) & 0.190 & 100 & 0.192 & 100 \\
\hline \hline
\multirow{3}{0.5cm}{0.5} & 50 & (0.357;0.643) & 0.442 & 100 & 0.471 & 100 \\
\cline{2-7}
 & 100 & (0.332;0.668) & 0.461 & 100 & 0.475 & 100 \\
\cline{2-7}
 & 500 & (0.334;0.666) & 0.486 & 100 & 0.489 & 100 \\
\hline \hline
\multirow{3}{0.5cm}{0.8} & 50 & (0.364;0.636) & 0.711 & 99 & 0.744 & 99 \\
\cline{2-7}
 & 100 & (0.352;0.648) & 0.746 & 99 & 0.764 & 99 \\
\cline{2-7}
 & 500 & (0.343;0.657) & 0.790 & 100 & 0.793 & 100 \\
\hline \hline
\multirow{3}{0.5cm}{0.9} & 50 & (0.438;0.562) & 0.795 & 93 & 0.846 & 94 \\
\cline{2-7}
 & 100 & (0.410;0.590) & 0.849 & 100 & 0.870 & 100 \\
\cline{2-7}
 & 500 & (0.336;0664) & 0.893 & 100 & 0.896 & 100 \\
\hline
\end{tabular}
\caption{Results obtained with $\pi=(\tfrac{1}{3},\tfrac{2}{3})$}
\label{tab:sim2}
\end{center}
\end{table}

\begin{table}
\begin{center}
\begin{tabular}{|c|c||c||c|c||c|c|}
\hline
$\alpha$ & $n$ & $\widehat{\pi}$ & $\widehat{\alpha_1}$ & $m_1$ & $\widehat{\alpha_2}$ & $m_2$ \\
\hline \hline
\multirow{3}{0.5cm}{0.1} & 50 & (0.324;0.345;0.331) & 0.120 & 74 & 0.132 & 75\\
\cline{2-7}
 & 100 & (0.337;0.339;0.324) & 0.099 & 82 & 0.103 & 84 \\
\cline{2-7}
 & 500 & (0.333;0.329;0.338) & 0.104 & 100 & 0.105 & 100 \\
\hline \hline
\multirow{3}{0.5cm}{0.2} & 50 & (0.345;0.334;0.321) & 0.176 & 95 & 0.184 & 98\\
\cline{2-7}
 & 100 & (0.327;0.342;0.331) & 0.186 & 99 & 0.193 & 99 \\
\cline{2-7}
 & 500 & (0.330;0.336;0.334) & 0.199 & 100 & 0.201 & 100 \\
\hline \hline
\multirow{3}{0.5cm}{0.5}&50 & (0.328;0.315;0.357) & 0.443 & 100 & 0.447 & 100\\
\cline{2-7}
 & 100 & (0.345;0.336;0.319) & 0.477 & 100 & 0.483 & 100 \\
\cline{2-7}
 & 500 & (0.340;0.323;0.332) & 0.495 & 100 & 0.496 & 100 \\
\hline \hline
\multirow{3}{0.5cm}{0.8}&50& (0.366;0.309;0.325) & 0.737 & 100 & 0.713 & 100\\
\cline{2-7}
 & 100 & (0.382;0.309;0.309) & 0.756 & 100 & 0.758 & 100 \\
\cline{2-7}
 & 500 & (0.336;0.328;0.336) & 0.793 & 100 & 0.795 & 100 \\
\hline \hline
\multirow{3}{0.5cm}{0.9}&50& (0.463;0.254;0.283) & 0.776 & 99 & 0.722 & 100\\
\cline{2-7}
 & 100 & (0.411;0.283;0.306) & 0.860 & 100 & 0.837 & 100 \\
\cline{2-7}
 & 500 & (0.348;0.322;0.330) & 0.895 & 100 & 0.893 & 100 \\
\hline
\end{tabular}
\caption{Results obtained with $\pi=(\tfrac{1}{3},\tfrac{1}{3},\tfrac{1}{3})$}
\label{tab:sim3}
\end{center}
\end{table}

\begin{table}
\begin{center}
\begin{tabular}{|c|c||c||c|c||c|c|}
\hline
$\alpha$ & $n$ & $\widehat{\pi}$ & $\widehat{\alpha_1}$ & $m_1$ & $\widehat{\alpha_2}$ & $m_2$ \\
\hline \hline
\multirow{3}{0.5cm}{0.1}&50 & (0.255;0.492;0.253) & 0.114 & 82 & 0.128 & 82 \\
\cline{2-7}
 & 100 & (0.248;0.495;0.257) & 0.097 & 89 & 0.102 & 90 \\
\cline{2-7}
 & 500 & (0.249;0.501;0.250) & 0.098 & 100 & 0.099 & 100 \\
\hline \hline
\multirow{3}{0.5cm}{0.2}&50 & (0.265;0.484;0.251) & 0.174 & 92 & 0.179 & 94 \\
\cline{2-7}
 & 100 & (0.255;0.494;0.251) & 0.177 & 99 & 0.182 & 99 \\
\cline{2-7}
 & 500 & (0.246;0.506;0.248) & 0.194 & 100 & 0.195 & 100 \\
\hline \hline
\multirow{3}{0.5cm}{0.5}&50& (0.279;0.467;0.254) & 0.467 &  100 & 0.467 & 100\\
\cline{2-7}
 & 100 & (0.250;0.494;0.256) & 0.466 & 100 & 0.467 & 100 \\
\cline{2-7}
 & 500 & (0.253;0.497;0.250) & 0.493 & 100 & 0.493 & 100 \\
\hline \hline
\multirow{3}{0.5cm}{0.8}&50& (0.299;0.496;0.205) & 0.716 & 100 & 0.693 & 100\\
\cline{2-7}
 & 100 & (0.282;0.472;0.246) & 0.766 & 100 & 0.765 & 100 \\
\cline{2-7}
 & 500 & (0.258;0.503;0.239) & 0.799 & 100 & 0.799 & 100 \\
\hline \hline
\multirow{3}{0.5cm}{0.9}&50& (0.356;0.422;0.222) & 0.802 & 100 & 0.715 & 100\\
\cline{2-7}
 & 100 & (0.332;0.460;0.208) & 0.858 & 100 & 0.827 & 100 \\
\cline{2-7}
 & 500 & (0.248;0.493;0.259) & 0.893 & 100 & 0.892 & 100 \\
\hline
\end{tabular}
\caption{Results obtained with $\pi=(\tfrac{1}{4},\tfrac{1}{2},\tfrac{1}{4})$}
\label{tab:sim4}
\end{center}
\end{table}

\newpage

\begin{table}
\begin{center}
\begin{tabular}{|c||c||c|c||c|c||c|c||c|c|}
\hline
 & & \multicolumn{2}{c||}{Ayes 2200} & \multicolumn{2}{c||}{Montgen\`evre 2200}  & \multicolumn{2}{c||}{N\'evache 1800} & \multicolumn{2}{c|}{Prorel 1800} \\
\hline
\hline
Model & Nb param & AIC & Nb obs & AIC & Nb obs & AIC & Nb obs & AIC & Nb obs \\
\hline
Indep. & 6 & {\bf 119.63} & 24 & {\bf 143.25} & 29 & {\bf 143.36} & 27 & {\bf 128.43} & 24 \\
\hline
Lag 1 & 42 & 141.26 & 22 & 157.58 & 27 & 169.06 & 22 & 154.02 & 20 \\
\hline
Lags 1-2 & 78 & 180.27 & 20 & 187.00 & 25 & 212.72 & 17 & 179.23 & 17  \\
\hline
\end{tabular}
\caption{Categorical time-series regression models applied to annual larch cones production}
\label{tab:larch-reg1}
\end{center}
\end{table}

\begin{table}
\begin{center}
\begin{tabular}{|c||c||c|c||c|c||c|c||c|c|}
\hline
 & & \multicolumn{2}{c||}{Ayes 2200} & \multicolumn{2}{c||}{Montgen\`evre 2200}  & \multicolumn{2}{c||}{N\'evache 1800} & \multicolumn{2}{c|}{Prorel 1800} \\
\hline
\hline
Model & Nb param & AIC & Nb obs & AIC & Nb obs & AIC & Nb obs & AIC & Nb obs \\
\hline
Indep. & 6 & 119.63 & 24 & 143.25 & 29 & 143.36 & 27 & {\bf 128.43} & 24  \\
\hline
Lag 1 & 12 & 116.34 & 22 & 136.3 & 27 & {\bf 139.32} & 22 & 136.62 & 20 \\
\hline
Lags 1-2 & 18 & {\bf 114.3} & 20 & {\bf 126.92} & 25 & 150.13 & 17 & 144.82 & 17 \\
\hline
\end{tabular}
\caption{Ordinal time-series regression models applied to annual larch cones production}
\label{tab:larch-reg2}
\end{center}
\end{table}

\begin{table}
\begin{center}
\begin{tabular}{|c||c|c|c||c|}
\hline
Valley & $\widehat{\pi}$ & $\widehat{\alpha_1}$ & $\widehat{\beta}$ & AIC \\
\hline \hline
Ayes 2200 & (0.167;0.042;0.292;0.292;0.125;0.083) & 0.082 & 0.774 & {\bf 121.06 } \\
\hline
Montgen\`evre 2200 & (0.138;0.069;0.172;0.310;0.241;0.069) & 0.070 & 0.935 & {\bf 118.70} \\
\hline
N\'evache 1800 & (0.185;0.185;0.074;0.185;0.296;0.074) & 0.032 & 0.871 & {\bf 125.94} \\
\hline
Prorel 1800 & (0.292;0.167;0.125;0.208;0.167;0.042) & 0.161 & 0.774 & {\bf 122.89} \\
\hline
\end{tabular}
\caption{DAR models applied to annual larch cones production}
\label{tab:larch-dar}
\end{center}
\end{table}

\newpage

\begin{table}
\begin{center}
\begin{tabular}{|c||c||c|c||c|c||c|c||c|c||c|c|}
\hline
 & & \multicolumn{2}{c||}{1987} & \multicolumn{2}{c||}{1988}  & \multicolumn{2}{c||}{1989} & \multicolumn{2}{c||}{1990}  & \multicolumn{2}{c|}{1987-1990} \\
\hline
\hline
Model & Nb param & AIC & Nb obs & AIC & Nb obs & AIC & Nb obs & AIC & Nb obs & AIC & Nb obs \\
\hline
Indep. & 3 & 100.56 & 47 & 111.52 & 45 & 121.61 & 48 & 133.6 & 48 & 468.21 &  188 \\
\hline
Lag 1 & 12 & {\bf 73.04} & 43 & {\bf 89.83} & 42 & {\bf 96.06} & 45 & {\bf 99.69} & 46 & 331.81 &  177 \\
\hline
Lags 1-2 & 21 & 79.52 & 40 & 93.57 & 37 & 97.17 & 42 & 102.90 & 44 & {\bf 315.96} & 167\\
\hline
\end{tabular}
\caption{Categorical time-series regression models applied to weekly planktonic abundance}
\label{tab:menard-reg1}
\end{center}
\end{table}

\begin{table}
\begin{center}
\begin{tabular}{|c||c||c|c||c|c||c|c||c|c||c|c|}
\hline
 & & \multicolumn{2}{c||}{1987} & \multicolumn{2}{c||}{1988}  & \multicolumn{2}{c||}{1989} & \multicolumn{2}{c||}{1990} & \multicolumn{2}{c|}{1987-1990} \\
\hline
\hline
Model & Nb param & AIC & Nb obs & AIC & Nb obs & AIC & Nb obs & AIC & Nb obs & AIC & Nb obs \\
\hline
Indep. & 3 & 100.56 & 47 & 111.52 & 45 & 121.61 & 48 & 133.6 & 48 & 468.21 & 188\\
\hline
Lag 1 & 6 & 65.77 & 43 & {\bf 84.78} & 42 & 95.67 & 45 & {\bf 92.66} & 46 & 323.5 &  177\\
\hline
Lags 1-2 & 9 & {\bf 63.69} & 40 & NA & 37 & {\bf 85.92} & 42 & 93.65 & 44 & {\bf 298.17} & 167\\
\hline
\end{tabular}
\caption{Ordinal time-series regression models applied to weekly planktonic abundance}
\label{tab:menard-reg2}
\end{center}
\end{table}

\begin{table}
\begin{center}
\begin{tabular}{|c||c|c|c||c|}
\hline
Year & $\widehat{\pi}$ & $\widehat{\alpha_1}$ & $\widehat{\beta}$ & AIC \\
\hline \hline
1987 & (0.625;0.167;0.042;0.000;0.167) & 0.540 & 0.923 & {\bf 108.78} \\
\hline
1988 & (0.489;0.311;0.089;0.000;0.111) & 0.308 & 0.865 & {\bf 143.11} \\
\hline
1989 & (0.354;0.417;0.167;0.000;0.062) & 0.445 & 0.923 & {\bf 132.68} \\
\hline
1990 & (0.375;0.271;0.146;0.000;0.208) & 0.484 & 0.923 & {\bf 135.10} \\
\hline
\hline
1987-1990 & (0.463;0.293;0.112;0.000;0.133) & 0.468 & 0.904 & {\bf 511.00} \\
\hline
\end{tabular}
\caption{DAR models applied to weekly planktonic abundance}
\label{tab:menard-dar}
\end{center}
\end{table}

\end{document}